\documentclass[12pt]{amsart}
\usepackage{lmodern}
\usepackage{amsthm}
\usepackage{mathptmx}
\usepackage{amsmath}
\usepackage{amssymb}
\usepackage{mathrsfs}
\usepackage{caption}
\usepackage[svgnames]{xcolor}
\usepackage{tikz}
  \usetikzlibrary{patterns}
\usepackage{graphicx,adjustbox}
\usepackage[letterpaper, portrait, margin=1in]{geometry}
\usepackage{csquotes}
\usepackage{amsthm}
\usepackage{mathtools}
\usepackage{wrapfig}
\usepackage{float}
\usepackage{subcaption}
\usepackage{amsrefs}

\MakeOuterQuote{"}

\newcommand{\fakecoprod}{\mathbin{\rotatebox[origin=c]{180}{$\Pi$}}}
\newcommand{\underln}[1]{\underline{\smash{#1}}}
\newcommand{\pa}{\partial}
\newcommand{\pr}{\prime}
\newcommand{\ra}{\rightarrow}
\newcommand{\lra}{\looparrowright}
\newcommand{\al}{\alpha}
\newcommand{\bb}{\beta}
\newcommand{\ga}{\gamma}
\newcommand{\de}{\delta}
\newcommand{\ep}{\epsilon}
\newcommand{\si}{\sigma}
\newcommand{\te}{\theta}

\newtheorem{thm}{Theorem}

\begin{document}

\title{Dehn's Lemma for Immersed Loops}

\author{Michael Freedman}
\address{\hskip-\parindent
        Michael Freedman\\
        Microsoft Research, Station Q, and Department of Mathematics\\
        University of California, Santa Barbara\\
        Santa Barbara, CA 93106\\}
\email{mfreedman@math.ucsb.edu}

\author{Martin Scharlemann}
\address{\hskip-\parindent
        Martin Scharlesmann\\
        Department of Mathematics\\
        University of California, Santa Barbara\\
        Santa Barbara, CA 93106\\}
\email{mgscharl@math.ucsb.edu}

\maketitle
\thispagestyle{myheadings}
\vspace{-1em}
We know five theorems whose conclusion is the existence of an embedded disk, perhaps with additional structure, in some larger space. Each introduced an influential technique and had broad consequences. They are: (1) - 1913 Boundary continuity of the Riemann mapping to Jordan domains (Carath\'{e}odory~\cite{Cara13}, Osgood and Taylor~\cite{Osgood13}); Carath\'{e}odory's proof introduced "external length"; applications to quasifuschian groups. (2) - 1944, the Whitney disk (with appropriate normal frame extension)~\cite{Whitney92}; applications to Whitney embedding theorem, h-cobordism theorem. (3) - 1957 Dehn's Lemma - loop theorem (Papakyriakopoulos~\cite{Papa57}); correctly treated triple points; applications: hierarchy in Haken manifolds, Thurston's geometrization theorem. (4) - 1982 Disk embedding theorem (Freedman~\cite{Freedman82}); used decomposition theory to identify Casson handles; application: topological classification of simply-connected 4-manifolds. (5) - Existence of Pseudoholomorphic disks (Gromov~\cite{Gromov85}); brought Kahler manifold techniques into symplectic context; applications: the nonsqueezing theorem, Seiberg-Witten invariants, quantum cohomology. This paper is a comment on (3); we prove a simply stated extension of Dehn's Lemma.

Let $M$ be a 3-manifold with boundary $\pa M$, and $\de: S^1 \lra \pa M$ be a generic immersion, where $S^1 \coloneqq [0, 2\pi] /$\tiny$0 \equiv 2\pi$\normalsize. By generic, we mean that $\de$ has only simple crossings. We say that $\de^\pr$ is $\de$ "displaced by a height function $f$" if $\de^\pr(\te) = (\de(\te), f(\te))$, where $f: S^1 \ra (0, \ep)$ is a Morse function and $[0, \ep]$ is a normal collar coordinate on $\pa M$ into $M$. We call a simple closed curve (scc) $\al \subset \textrm{int} (M)$ unknotted iff $\al$ bounds an embedded disk $\bar{\al}: D^2 \ra \textrm{int} (M)$.

\begin{thm}
  Let $\de: S^1 \lra \pa M$ be a generic immersed loop so that the composition into $M$ is null homotopic. There is a height function $f: S^1 \ra (0, \ep)$ so that $\de^\pr = (\de, f): S^1 \hookrightarrow \textnormal{int} (M)$ is unknotted.
\end{thm}

The theorem readily implies two familiar facts:

\begin{enumerate}
  \item Dehn's Lemma: For any scc $\de \subset \pa M$ which is null homotopic in $M$ there is a properly embedded disk $(D, \pa D) \subset (M, \pa M)$ with $\pa D$ parameterizing $\de$.
  \vspace{0.5em}
  \item Given any knot diagram there is always a way to rechoose the crossings to produce an unknot. (This is the case when $(M, \pa M) \cong (B^3, \pa B^3)$ is a 3-ball.)
\end{enumerate}

\underln{Regarding 1}: A short argument connects the special case of the theorem where $\de$ is one-to-one to Dehn's Lemma. Let $C$ be the annular collar joining $\de$ to $\de^\pr$ and $D$ a disk with boundary $\de^\pr$. If $C$ and $D$ are in general position initially there may be arcs of intersection, but a perturbation of $D$ near $\pa D = \de^\pr$ starting from an outermost arc ensures that $C$ and $D$ intersect only in sccs contained in int$(C)$ and int$(D)$. Let $\si \subset D$ be an innermost circle of intersection bounding a subdisk $\overline{\si} \in D$. $\si$ may be paired with either an essential or inessential scc in $C$. Perform disk exchanges to modify $C$ until either $C \cup D$ is an embedded disk or some innermost $\si$ is paired to an essential scc in $C$. In this case a final cut and glue operation yields an embedded disk with boundary $\de$.

\underln{Regarding 2}: We should note a subtlety. The height function produced in the theorem may, in general, be more complicated than the familiar height function which solves the knot diagram problem (see Figure 1b). For knot diagrams the unknotting function may be taken to be any function with a unique local maximum and unique local minimum.

Actually Theorem 1 is a corollary of a stronger Theorem 0, better adapted to the required induction.

\noindent
\textbf{Theorem 0.} \textit{Let $\de: S^1 \lra \pa M$ be a generic immersed loop with base point $\ast$ ($\ast$ is assumed disjoint from multiple points of $\de$) whose composition into $M$ is null homotopic. There is a height function $f: S^1 \ra [0, \epsilon)$, $f(\ast) = 0$, $f(S^1 \backslash \ast) \subset (0, \epsilon)$, so that $\de^\pr = \{ (\de (\te), f(\te)\}$ bounds an embedded disk $\Delta \subset M$ with $\Delta \cap \pa M = \ast$. We say $\de^\pr$ collapses to $\ast$, and call $\Delta$ a "lollypop" for $\de^\pr$.}

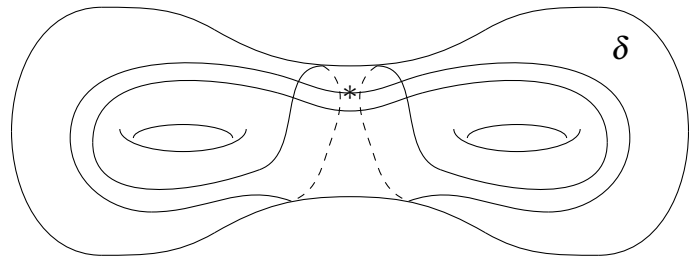
\begin{figure}[ht]
  \centering
  \begin{subfigure}{\textwidth}
  \centering
  \begin{tikzpicture}[scale = 1.2]
    \draw (0,0.1) to [out = 0, in = 210] (1.5,0.5);
    \draw (1.5,0.5) to [out = 30, in = 180] (2.75,0.75);
    \draw (0,-1.35) to [out=0, in = 150] (1.5,-1.75);
    \draw (1.5,-1.75) to [out=-30, in=180] (2.75,-2);
    \draw (2.75,0.75) to [out=0, in = 90] (3.75, -0.625);
    \draw (3.75,-0.625) to [out=-90, in=0] (2.75, -2);

    \draw (0,0.1) to [out=180,in=-30] (-1.5,0.5);
    \draw (-1.5,0.5) to [out=150, in=0] (-2.75, 0.75);
    \draw (0,-1.35) to [out=180,in=30] (-1.5,-1.75);
    \draw (-1.5,-1.75) to [out=210,in=0] (-2.75,-2);
    \draw (-2.75,0.75) to [out=180,in=90] (-3.75,-0.625);
    \draw (-3.75,-0.625) to [out=-90,in=180] (-2.75, -2);

    \draw (1.15,-0.6) arc (-180:0:0.7 and 0.25); 
    \draw (-1.15,-0.6) arc (0:-180:0.7 and 0.25);
    \draw (1.3, -0.7) arc (180:0:0.55 and 0.15); 
    \draw (-1.3,-0.7) arc (0:180:0.55 and 0.15);

    \draw (0.3,0.1) to [out=0, in=160] (1, -1.05);
    \draw (1,-1.05) to [out=-20, in=-90] (2.85,-0.75);
    \draw (2.85,-0.75) to [out=90,in=15] (0.8,-0.2);
    \draw (0.8,-0.2) to [out=195, in=0] (0,-0.4);
    \draw (0.3,0.1) [dashed] to [out=210, in=110] (0.2,-0.7);
    \draw (0.2,-0.7) [dashed] to [out=-70, in=160] (0.65,-1.4);
    \draw (0.65,-1.4) to [out=20,in=180] (2.2,-1.52);
    \draw (2.2,-1.52) to [out=0,in=-90] (3.1,-0.7);
    \draw (3.1,-0.7) to [out=90,in=15] (0.8,0);
    \draw (0.8,0) to [out=195,in=0] (0,-0.2);

    \draw (-0.3,0.1) to [out=180, in=20] (-1, -1.05);
    \draw (-1,-1.05) to [out=200, in=-90] (-2.85,-0.75);
    \draw (-2.85,-0.75) to [out=90,in=165] (-0.8,-0.2);
    \draw (-0.8,-0.2) to [out=-15, in=180] (-0,-0.4);
    \draw (-0.3,0.1) [dashed] to [out=-30, in=70] (-0.2,-0.7);
    \draw (-0.2,-0.7) [dashed] to [out=250, in=20] (-0.65,-1.4);
    \draw (-0.65, -1.4) to [out=160, in=0] (-2.2, -1.52);
    \draw (-2.2,-1.52) to [out=180,in=-90] (-3.1,-0.7);
    \draw (-3.1,-0.7) to [out=90,in=165] (-0.8,0);
    \draw (-0.8,0) to [out=-15,in=180] (0,-0.2);
    \node at (0,-0.2) {$\ast$};
    \node at (3,0.3) {$\delta$};
  \end{tikzpicture}
  \caption{This example of a loop on the boundary of a genus 2 handlebody shows that the height function $f$, constructed in the proof of Theorem 0, cannot always have a unique local maximum, the familiar form for unknotting knot diagrams on the 2-sphere. The unknotting function $f$, relative to the indicated base point $\ast$ with $f(\ast) = 0$, must have at least two local maxima.}
  \end{subfigure}\hfill
  
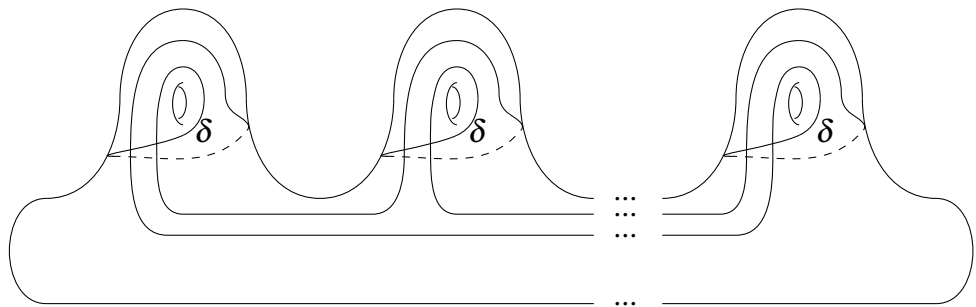
\vspace{1em}
  \begin{subfigure}{\textwidth}
  \centering
  \begin{tikzpicture}[scale = 0.7]
    \draw (0,0) to [out=180,in=90] (-1.2, -1.7);
    \draw (-1.2, -1.7) to [out=-90, in=0] (-2.6,-3.6);
    \draw (0,0) to [out=0, in=90] (1.2, -1.7);
    \draw (1.2, -1.7) to [out=-90, in=180] (2.6,-3.6);
    \draw (0, -1.4) arc (90:270:0.2 and 0.4);
    \draw (-0.09, -1.48) arc (90:-90:0.15 and 0.3);
    \draw (0, -1.1) to [out=180, in=75] (-0.4, -1.5);
    \draw (-0.4, -1.5) to [out=255, in=90] (-0.5, -2.8);
    \draw (-1, -2.8) to [out=90, in=180] (0,-0.6);
    \draw (0, -1.1) to [out=0, in=90] (0.4, -1.7);
    \draw (0.4, -1.7) to [out=-90, in=25] (0, -2.4);
    \draw (0, -2.4) to [out=205, in=180] (-1.32, -2.8);
    \draw [dashed] (-1.34, -2.8) to [out=0, in=230] (1.22,-2.3);
    \draw (0, -0.6) to [out=0, in=90] (0.8,-1.6);
    \draw (0.8,-1.6) to [out=-90, in=145] (1,-2);
    \draw (1, -2) to [out=-35, in=50] (1.22, -2.3);

    \draw (-5.2,0) to [out=180,in=90] (-6.4,-1.7);
    \draw (-6.4,-1.7) to [out=-90, in=0] (-7.8,-3.6);
    \draw (-5.2,0) to [out=0, in=90] (-4,-1.7);
    \draw (-4,-1.7) to [out=-90, in=180] (-2.6,-3.6);
    \draw (-5.2,-1.4) arc (90:270:0.2 and 0.4);
    \draw (-5.29,-1.48) arc (90:-90:0.15 and 0.3);
    \draw (-5.2,-1.1) to [out=180, in=75] (-5.6,-1.5);
    \draw (-5.6,-1.5) to [out=255, in=90] (-5.7,-2.8);
    \draw (-6.2,-2.8) to [out=90, in=180] (-5.2,-0.6);
    \draw (-5.2,-1.1) to [out=0, in=90] (-4.8,-1.7);
    \draw (-4.8,-1.7) to [out=-90, in=25] (-5.2,-2.4);
    \draw (-5.2,-2.4) to [out=205, in=180] (-6.52,-2.8);
    \draw [dashed] (-6.54,-2.8) to [out=0, in=230] (-3.98,-2.3);
    \draw (-5.2,-0.6) to [out=0, in=90] (-4.4,-1.6);
    \draw (-4.4,-1.6) to [out=-90, in=145] (-4.2,-2);
    \draw (-4.2,-2) to [out=-35, in=50] (-3.98,-2.3);

    \draw (-5.7, -2.8) to [out=-90, in=180] (-5.2,-3.9) to (-1.6, -3.9) to [out=0, in=-90] (-1, -2.8);
    \draw (-0.5, -2.8) to [out=-90, in=180] (0, -3.9) to (2.6, -3.9);
    \draw (-6.2, -2.8) to [out=-90, in=180] (-5.5, -4.3) to (2.6, -4.3);
    \draw (-7.8, -3.6) to [out=180, in=90] (-8.5, -4.6) to [out=-90, in=180] (-7.8, -5.6) -- (2.6, -5.6);

    \draw (6.5,0) to [out=180,in=90] (5.3,-1.7);
    \draw (5.3,-1.7) to [out=-90, in=0] (3.9,-3.6);
    \draw (6.5,0) to [out=0, in=90] (7.7,-1.7);
    \draw (7.7,-1.7) to [out=-90, in=180] (9.1,-3.6);
    \draw (6.5,-1.4) arc (90:270:0.2 and 0.4);
    \draw (6.41,-1.48) arc (90:-90:0.15 and 0.3);
    \draw (6.5,-1.1) to [out=180, in=75] (6.1,-1.5);
    \draw (6.1,-1.5) to [out=255, in=90] (6,-2.8);
    \draw (5.5,-2.8) to [out=90, in=180] (6.5,-0.6);
    \draw (6.5,-1.1) to [out=0, in=90] (6.9,-1.7);
    \draw (6.9,-1.7) to [out=-90, in=25] (6.5,-2.4);
    \draw (6.5,-2.4) to [out=205, in=180] (5.18,-2.8);
    \draw [dashed] (5.17,-2.8) to [out=0, in=230] (7.72,-2.3);
    \draw (6.5,-0.6) to [out=0, in=90] (7.3,-1.6);
    \draw (7.3,-1.6) to [out=-90, in=145] (7.5,-2);
    \draw (7.5,-2) to [out=-35, in=50] (7.72,-2.3);

    \node at (3.2,-3.6) {...};
    \node at (3.2,-3.9) {...};
    \node at (3.2,-4.3) {...};
    \node at (3.2,-5.6) {...};
    \draw (5.5, -2.8) to [out=-90, in=0] (5, -3.9) -- (3.9, -3.9);
    \draw (6, -2.8) to [out=-90, in=0] (5.3, -4.3) -- (3.9, -4.3);
    \draw (3.9, -5.6) to (9.1, -5.6);
    \draw (9.1, -5.6) to [out=0, in=-90] (9.8, -4.6) to [out=90, in=0] (9.1, -3.6);
    \node at (0,-6.6) {handlebody genus $= m+1$};
    \node at (-4.8,-2.3) {$\delta$};
    \node at (0.4,-2.3) {$\delta$};
    \node at (7,-2.3) {$\delta$};
    \end{tikzpicture}
    \caption{For this absolute example, the unknotting function $f$ produced by Theorem 1 must have at least $m$ local minima (and maxima)}
    \end{subfigure}
    \caption{}
\end{figure}

\begin{proof}
  No methods post dating Papakyriakopoulos \cite{Papa57} are required. He would have found this proof rather easy to understand and perhaps to generate. First we build a tower.

  $\de$ is the pointed, immersed loop, $D$ the immersed null homotopy bounding it, and $N$ the regular neighborhood of $D$ in $M$. Subscript will indicate height in the tower.

  To build the tower we should ask about the $Z_2$ Betti number $b_1(N) \coloneqq \textrm{rank} (H_1 (N; Z_2))$. If $b_1(N) = 0$, there is no tower. For homological reasons $\pa N$ is a disjoint union of 2-spheres and $\de$ is contained in one of these: $\de \lra S \subset \pa N$.

  If $b_1(N) > 0$ choose a 2-fold cover $\widetilde{N}_0 \ra N_0 \coloneqq N$ and choose a lift $l_1: D \lra \widetilde{N}_0$ and let $N_1 = \textrm{neib} (l_1(D)) \subset \widetilde{N}_0$. If $b_1(N_1) = 0$ then $N_1$ is the top of the tower. If $b_1(N_1) > 0$ continue and find a 2-fold cover $\widetilde{N}_1 \ra N_1$, and lift $l_2: D \lra \widetilde{N}_1$ and set $N_2 = \textrm{neib} (l_2(D)) \subset \widetilde{N}_1$. Again if $b_1(N_2) = 0$ then $N_2$ is the top of the tower. If not, proceed to construct $l_3 ... l_n$ and $N_3 ... N_n$ until $b_1(N_n) = 0$. A simple complexity argument, where the complexity can be the number of simplicies identified by $l_k$ for a fixed triangulation of $D$ making all $l_k$ simplicial, shows that the tower is indeed finite.

  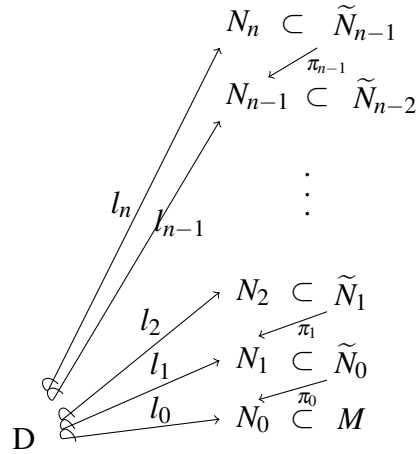
\begin{figure}[ht]
    \centering
    \begin{tikzpicture}[scale = 1.3]
      \node at (0,0) {D};
      \draw [->] (0.25,0.5) to (2,4);
      \node at (2.25, 4.25) {$N_n$};
      \node at (2.75, 4.25) {$\subset$};
      \node at (3.5, 4.25) {$\widetilde{N}_{n-1}$};

      \draw [->] (0.3, 0.4) to (2, 3.25);
      \node at (2.4, 3.5) {$N_{n-1}$};
      \node at (3, 3.5) {$\subset$};
      \node at (3.7, 3.5) {$\widetilde{N}_{n-2}$};
      \draw [->] (3,4) to (2.5, 3.7);
      \node at (2.9, 2.7) {.};
      \node at (2.9, 2.5) {.};
      \node at (2.9, 2.3) {.};

      \draw [->] (0.4, 0.2) to (2, 1.5);
      \node at (2.35, 1.5) {$N_2$};
      \node at (2.85, 1.5) {$\subset$};
      \node at (3.35, 1.5) {$\widetilde{N}_{1}$};
      \draw [->] (3.1, 1.3) to (2.4, 1.05);

      \draw [->] (0.4, 0.1) to (2, 0.8);
      \node at (2.35, 0.8) {$N_1$};
      \node at (2.85, 0.8) {$\subset$};
      \node at (3.35, 0.8) {$\widetilde{N}_{0}$};
      \draw [->] (3.1, 0.6) to (2.4, 0.4);

      \draw [->] (0.4, 0) to (2, 0.2);
      \node at (2.35, 0.2) {$N_0$};
      \node at (2.85, 0.2) {$\subset$};
      \node at (3.35, 0.2) {$M$};

      \node at (1, 2.4) {$l_n$};
      \node at (1.6, 2.2) {$l_{n-1}$};
      \node at (1.3, 1.2) {$l_2$};
      \node at (1.4, 0.75) {$l_1$};
      \node at (1.4, 0.3) {$l_0$};

      \draw (0.35, 0.56) [out = 135, in = 45] to (0.2, 0.6);
      \draw (0.2, 0.6) [out = 225, in = 225] to (0.25, 0.5);

      \draw (0.4, 0.46) [out = 135, in = 45] to (0.25, 0.5);
      \draw (0.25, 0.5) [out = 225, in = 225] to (0.3, 0.4);

      \draw (0.52, 0.22) [out = 135, in = 45] to (0.37, 0.3);
      \draw (0.37, 0.3) [out = 225, in = 225] to (0.4, 0.2);

      \draw (0.52, 0.10) [out = 120, in = 60] to (0.39, 0.18);
      \draw (0.39, 0.18) [out = 210, in = 210] to (0.4, 0.1);

      \draw (0.53, -0.03) [out = 90, in = 60] to (0.39, 0.08);
      \draw (0.39, 0.08) [out = 190, in = 190] to (0.4, 0);

      \node at (3.1,3.8) {\tiny $\pi_{n-1}$};
      \node at (2.9,1.1) {\tiny $\pi_{1}$};
      \node at (2.9,0.4) {\tiny $\pi_{0}$};
  \end{tikzpicture}
  \caption{The tower}
  \end{figure}

  \underln{Observation:} $\de_n \coloneqq l_n(\pa D)$ can be unknotted by a suitable resolution of its crossings (i.e. a normal function $f$ as in Theorem 1). As in Theorem 0, given for any base point $\ast$ (chosen away from crossings) we can resolve crossings $\de_n \ra \de_n^\pr$ and produce an embedded disk $\Delta$ with $\pa \Delta = \de_n^\pr$ and $\Delta \cap \pa N_n = \ast$.

   \textit{Explanation of $\de_n^\pr$ and its null-isotopy.} One might expect to choose $\de_n^\pr = \{\de_n(\te), f(\te)\}$, $f(\ast)$ = 0, $f(\te \neq \ast) > 0$, where $f$ has a unique local maximum and a unique local minimum. However, this solution does \underln{not}, in general, push down the tower. We prefer to give a second solution. There are two cases. If $l_n$ is one-to-one $\de_n = l_n (\pa D)$ bounds a hemisphere $E_1 \subset S \subset \pa N_n$. Pushing the disk $E_1$ normally toward interior $N_n$ gives a disk $(\Delta, \ast) \subset (N_n, \pa N_n)$ bounding $\de_n^\pr$.

   Now assume $\de_n$ is not one-to-one. Let $\al \subset \de_n$, not containing $\ast$, be a subarc so that $l_n(\al)$ is a scc $\subset \pa N$. Let $E_1$ be one of the two disks in $\pa N_n$ bounded by $l_n (\al)$. Begin to resolve the crossings of $\de$ by following this rule: $l_n(\al)$ lies \underln{above} $l_n(\bb)$, where $\bb$ is the complementary arc, $\de = \al \cup \bb$ and $\al \cap \bb = \pa \al = \pa \bb$. Also the two endpoints of $\al$ are at a crossing: resolve this crossing arbitrarily.

   \begin{figure}[ht]
    \centering
    \begin{tikzpicture}[scale = 0.7]
      \draw [thick] (0,0) to [out=30, in=-90] (1, 1.5);
      \draw [thick] (1,1.5) to [out=90,in=0] (0,3);
      \draw [thick] (0,3) to [out=180, in=90] (-1,1.5);
      \draw [thick] (-1,1.5) to [out=-90,in=150] (0,0);

      \draw (0,0) to [out=-30, in=240] (3.3,1.1);
      \draw (3.3,1.1) to [out=60, in=45] (-1,3.6);
      \draw (-1,3.6) to [out=225, in=160] (1.7,-0.5);
      \draw (1.7,-0.5) to [out=-20, in=0] (1,-2.5);
      \draw (1,-2.5) to [out=180,in=270] (-1.3,-0.5);

      \draw [dotted, thick] (-1.3,-0.5) to [out=90,in=225] (-0.4,0.9);
      \draw [dotted, thick] (-0.4,0.9) to [out=45,in=-60] (0.2,2);
      \draw [dotted, thick] (0.2, 2) to [out=120, in=65] (-4, 1.6);
      \draw [dotted, thick] (-4,1.6) to [out=245, in=170] (-3.1,-0.6);
      \draw [dotted, thick] (-3.1, -0.6) to [out=-10, in=195] (-1.3,-0.5);
      \draw (-1.3,-0.5) to [out=15, in = 210] (0,0);

      \node at (1.5,2) {$\alpha$};
      \draw [->](1.5, 1.8) to [out=-90,in=0] (1.2, 1.7);

      \node at (4.2,1.5) {$\beta$};
      \draw [->] (4.05, 1.5) to (3.55, 1.5);
      \node at (-4.4, 0.9) {$\alpha^\prime$};
      \node at (2.7,3) {$\ast$};
      \node at (0.1, -3.1) {$\delta = \al \cup \beta$};
    \end{tikzpicture}
    \caption{}
  \end{figure}

   Now an isotopy across the disk $E_1$ bounding $\al$ effectively erases the loop $\al$. If the simplified diagram is a scc we may continue the unknotting using a hemisphere 2-cell $E_2 \subset S, \pa E_2 = l_n (\beta)$ making no further crossing choices. If the simplified diagram is still singular choose another arc $\al^\pr \subset \bb$ whose image is again a scc not containing $\ast$ and proceed as before. Continuing in this way, guided by a sequence of embedded 2-cells, say $\{E_1, .., E_j\}$, all crossings are eventually assigned so that the diagram resolution is unknotted with $\ast$ remaining in the final 2-cell $E_j$. The cells $E_1, ..., E_j$ determine a sequence of isotopies $\mathcal{I}_1, ..., \mathcal{I}_j$ so that the composition $\mathcal{I}_j \circ ... \circ \mathcal{I}_1$ shrinks $\de_n^\pr$ toward the base point $\ast$. Note: the 2-cells $E_i$ will not generally have disjoint interiors.
   \qed

  This solution will now be pushed down the tower

  \textit{Dissection of disks with double arcs}: Let $E$ be a properly embedded disk $(E, \pa E) \subset (M^3, \pa M)$ in a 3-manifold and $\pi: (M, \pa M) \ra (P, \pa P)$ a 2-fold cover with covering translation $t$. Assume $E$ and $t E$ are in generic position and meet only in double arcs and double loops. Double loops are easily removed by an innermost circle argument and will not be discussed further. $E$ is assumed to have a base point $\ast \subset \pa E$. Call $\pi(E) = F$ and $\pi(\ast) = \ast$ in an abuse of notation.

  We now describe how to form an ordered list of embedded disks $F_1, ..., F_k \subset (P, \pa)$ from pieces (some used several times) of $F$. Constructing $\{F_1, ..., F_k\}$ will dictate "crossing choices" for $\pa F \coloneqq \ga$, which yield $\ga^\pr$, bounding a lollypop. Each $F_i$ has base point $\ast_i$ with $\ast_k = \ast \subset \pa F$. Furthermore, once $\{F_1, ..., F_k\}$ are constructed we may view them as the instructions for an isotopy $\mathcal{I} = \mathcal{I}_k \circ ... \circ \mathcal{I}_1$, as above, collapsing $\ga^\pr$ toward $\ast$.

  Let $\bb_1$ be an outermost arc cutting off an outermost disk $D_1 \subset F$ not containing $\ast$. Let $\bb_1^\pr$ be the partner arc of $\bb_1$ and $D_1^\pr$ be the subdisk of $F$ cut off by $\bb_1^\pr$ which also does not contain $\ast$. There are two cases, shown in Figure 4.

  \begin{figure}[ht]
    \centering
    \begin{tikzpicture}
      \draw (0,0) ellipse (4 and 1);
      \path [pattern=north east lines] (-4,0) to [out=90,in=210] (-3.65,0.4) to [out=30, in=192] (-1.8,0.9) to (-1.8,-0.9) to [out=168, in=-30] (-3.65, -0.4) to [out=150, in=-90] (-4,0);
      \draw (-1.8,0.9) to (-1.8,-0.9);

      \path [pattern=north west lines] (4,0) to [out=90,in=-30] (3.65,0.4) to [out=150, in=-12] (1.8,0.9) to (1.8,-0.9) to [out=12, in=210] (3.65, -0.4) to [out=30, in=270] (4,0);
      \draw (1.8,0.9) to (1.8,-0.9);

      \node at (0,1) {$\ast$};
      \node at (-1.5,0) {$\beta_1^\prime$};
      \node at (1.5,0) {$\beta_1$};
      \node at (-3.2,-1.2) {$D_1^\prime$};
      \node at (2.7,1.3) {$D_1$};
      \node at (-1.8,1.1) {over};
      \node at (-1.8,-1.1) {over};
      \draw [->] (-3.2, -1) to (-3.1,-0.7);
      \draw [->] (2.6,1.1) to (2.4,0.85);

      \draw (9.2,0) ellipse (3 and 1);
      \node at (6.2,0) {$\ast$};

      \path [pattern = north east lines] (7.9,0.9) to [out=10, in=180] (9.2,1) to [out=0, in=170] (10.5,0.9) to (10.5,-0.9) to [out=190, in=0] (9.2,-1) to [out=180,in=-10] (7.9,-0.9) to (7.9, 0.9);
      \draw (7.9,0.9) to (7.9,-0.9);
      \draw (10.5, 0.9) to (10.5,-0.9);
      \path [pattern = crosshatch] (10.5, 0.9) to [out=-10,in=139] (12.05,0.3) to [out= -45, in=90] (12.2,0) to [out=-90, in=45] (12.05, -0.3) to [out=221, in=10] (10.5,-0.9);
      \node at (7.6,0) {$\beta_1^\prime$};
      \node at (10.5,-1.2) {$\beta_1$};
      \node at (7.9,1.1) {over};
      \node at (7.9,-1.1) {over};
      \node at (9.1,-1.5) {$D_1^\prime$};
      \node at (11.5,1.2) {$D_1$};
      \draw [->] (9.1,-1.3) to (9.2, -1.05);
      \draw [->] (11.4,1.05) to (11.2, 0.8);
      \node at (0,-2) {case 1};
      \node at (9.2,-2) {case 2};
    \end{tikzpicture}
    \caption{}
  \end{figure}

  In both cases $D_1 \cup D_1^\pr$ can be perturbed into a proper map of a disk $F_0$. There is a residual general position proper map with only double arc singularities of a disk formed from $D_1 \cup (F \backslash D_1^\pr)$. In both cases, label this map $F_1$. Because $D_1$ is outermost, both $F_0$ and $F_1$ have fewer double arcs than $F$.

  We need to discuss crossing choices (resolutions) and base point choices. The end points on $\pa D_1$, i.e. its crossings with \underln{other} segments of $\pa F_1$, are deemed \underln{overcrossings} (in both cases). $F_0$ is provided a base point $\ast_0$ in its copy of $D_1$ and $F_1$ retains $\ast_1 \coloneqq \ast$ as its base point.

  Now unless $F_0 (F_1)$ is embedded find an outermost arc cutting off an outermost disk $D$, with $D^\prime$ disjoint from $\ast_0 (\ast_1)$ and sharing a double arc with $D$, \underln{not} containing $\ast_0 (\ast_1)$ and as above dissect:

  
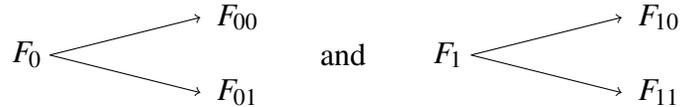
\begin{center}
  \begin{tikzpicture}
    \node at (0,0) {$F_0$};
    \draw [->] (0.3,0) to (2.3,0.5);
    \draw [->] (0.3,0) to (2.3,-0.5);
    \node at (2.8,0.5) {$F_{00}$};
    \node at (2.8,-0.5) {$F_{01}$};
    \node at (4.2,0) {and};
    \node at (5.6,0) {$F_1$};
    \draw [->] (5.9,0) to (7.9,0.5);
    \draw [->] (5.9,0) to (7.9,-0.5);
    \node at (8.4,0.5) {$F_{10}$};
    \node at (8.4,-0.5) {$F_{11}$};
  \end{tikzpicture}
  \end{center}

  \noindent
  to obtain general position proper pointed maps $F_{00}, F_{01}, F_{10},$ and $F_{11}$ of disks with only double arc singularities. We call such maps \underln{good} maps. In each case follow the preceding rule for resolving the crossing of $\pa D$ as \underln{overcrossings}.

  Continuing in this way a dyadic tree of good maps is obtained. The leaves of this tree are called \underln{great maps} as they have the additional property of being embeddings. (No disjointness has been constructed or assumed for these great maps.) The leaves are now linearly ordered by the base 2 numerical value of their subscripts considered as decimals. By the time we reach the leaves \underln{all} crossings have been resolved. These great disks, monotonically reindexed, become the ordered list $\{F_1, .., F_k\}$, built from pieces of $F_1$, which we sought. We call this list a \underln{great sequence} (for $\pa F$) \underln{guiding a collapse} of $(\pa F)^\pr$ to $\ast$. Note that the great sequence uniquely defines the crossing resolution of $\pa F^\pr$ of $\pa F$. The simple expedient of successively declaring $\pa D$ to over-cross other segments of the boundary has given us the well-defined crossing choices and thus defines $(\pa F)^\pr$. The isotopy $\mathcal{I}$ shrinks $(\pa F)^\pr$ toward $\ast$.

  \begin{figure}[H]
    \centering
    \begin{tikzpicture}
      \node at (0,0) {$F$};
      \draw [->] (0.2,0) node (v1) {} -- (1.7,0.5);
      \draw [->] (0.2,0) -- (1.7,-0.5);
      \node at (2,0.5) {$F_0$};
      \node at (2,-0.5) {$F_1$};
      \draw [->] (2.3,0.5) -- (3.8,0.5);
      \draw [->] (2.3,-0.5) -- (3.8,-0.5);
      \draw [->] (2.3, 0.7) -- (3.8,1.3);
      \draw [->] (2.3, -0.7) -- (3.8,-1.3);
      \node at (4.1,1.3) {$F_{00}$};
      \node at (4.1,0.5) {$F_{01}$};
      \node at (4.1,-0.5) {$F_{10}$};
      \node at (4.1,-1.3) {$F_{11}$};
      \draw [->] (4.5, 0.5) to (6,0.8);
      \draw [->] (4.5, 0.5) to (6,0.2);
      \draw [->] (4.5,-0.5) to (6,-0.2);
      \draw [->] (4.5,-0.5) to (6,-0.8);
      \node at (6.4,0.8) {$F_{010}$};
      \node at (6.4,0.2) {$F_{011}$};
      \node at (6.4,-0.2) {$F_{100}$};
      \node at (6.4,-0.8) {$F_{101}$};
      \draw (6.9,-0.2) -- (8.4,0);
      \draw (6.9,-0.2) -- (8.4,-0.4);
      \node at (8.9,0) {$F_{1000}$};
      \node at (8.9,-0.4) {$F_{1001}$};
    \end{tikzpicture}
    \caption{Sample tree with great maps as leaves}
  \end{figure}
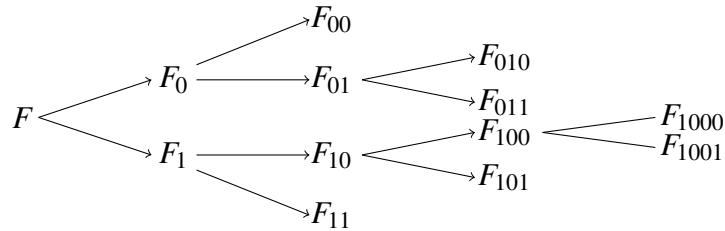

  \underln{To summarize the Dissection section}:

  The tree, as constructed above, of good disks terminating in great disks is a \underln{dissection} of $F$. The same tree, according to our convention, determines crossing choices $\ga^\pr = (\pa F)^\pr$ for $\ga = \pa F$. The great sequence $\{F_1, ..., F_k\}$ (the leaves of the dissection) are said to \underln{guide} a sequence of isotopes $\mathcal{I}_k \circ ... \circ \mathcal{I}_1 \coloneqq \mathcal{I}$ with $\mathcal{I}_i$ supported near $F_i$, $1 \leq i \leq k$, so that $\mathcal{I}$ collapses $\ga^\pr$ toward its base point.

  \noindent
  \textbf{Example:} We illustrate (Figure 6) dissection, the crossing resolutions, and the collapsing isotopy with an example. Our crossing convention implies that middle sheet of the three sheets of the outermost $D$ in $F_{...0} \fakecoprod F_{...1}$ should be compressed slightly into interior ($M$) near $\pa \bb_1$ when interpreting the disk $F_j$, $1 \leq j \leq k$, which contain this sheet of $D$ as guiding the isotopes $\mathcal{I}_j$.

  \begin{figure}
    \centering
    \begin{tikzpicture}[scale = 0.8]
    \draw (0,0) ellipse (4 and 2);
    \draw (-3.33, 1.11) to [out=-70, in=225] (-2.7,0.8);
    \draw (-2.7,0.8) to [out = 45, in=-70] (-2.7, 1.48);
    \draw (-2,1.73) -- (-2,-1.73);
    \draw (-1.3, 1.9) to [out=-45, in=90] (-0.3,-2);
    \draw (1.1,1.92) -- (1.1,-1.92);
    \draw (2.3,1.63) -- (2.3,-1.63);
    \draw (3, -1.3) to [out=135, in=225] (3,-0.9) to [out=45, in=150] (3.45, -1);
    \draw (-0.3,2) to [out=275, in=120] (0.05,0.4);
    \draw (0.45, 0.4) to [out=60, in=265] (0.7, 1.97);
    \draw (0.05, 0.4) to [out=-60, in=240] (0.45, 0.4);
    \draw (-0.1, 2) to [out=280, in=120] (0.1, 1.3);
    \draw (0.5, 1.99) to [out=260, in=60] (0.3, 1.3);
    \draw (0.1, 1.3) to [out=-60, in=240] (0.3, 1.3);

    \node at (-3.5,-1) {$\ast$};
    \node at (-2.3,1.2) {$d^\prime$};
    \node at (-2.2,-0.2) {$a^\prime$};
    \node at (-0.6,-0.2) {$b^\prime$};
    \node at (0.5,0.2) {$c^\prime$};
    \node at (0.35,2.37) {$d$};
    \draw [->] (0.3,2.2) -- (0.2,1.7);
    \node at (1.3,-0.1) {$b$};
    \node at (2.5,-0.1) {$a$};
    \node at (3.4,-0.7) {$c$};
    \node at (0.2,-2.5) {\large $F$};

    \draw (6.8, -2.1) to [out=60, in=-100] (8.4,0);
    \draw (8.4,0) to [out=70, in=-45] (7.6,1.5);
    \draw (7.6,1.5) to [out=135, in=-90] (6.9,2.4) to [out=90, in=180] (7.6,3.5);
    \draw (7.6,3.5) to [out=0, in=90] (8.4,2.4) to [out=-90, in=45] (7.6, 1.5);
    \draw (7.6, 1.5) to [out=225, in=110] (6.9,0);
    \draw (8.4, -2.1) to [out=120, in=280] (6.9,0);

    \draw [dashed](9.7,-1.4) to [out=60, in=-100] (11.2,0.7);
    \draw (11.2,0.7) to [out=70, in=-45] (10.5,2.2);
    \draw [dashed] (10.5,2.2) to [out=135, in=-90] (9.8,3.1) to [out=90, in=180] (10.5,4.2);
    \draw (10.5,4.2) to [out=0, in=90] (11.3,3.1) to [out=-90, in=45] (10.5,2.2);
    \draw [dashed] (10.5,2.2) to [out=225, in=110] (9.8,0.7);
    \draw [dashed](11.3,-1.4) to [out=120, in=280] (9.8,0.7);
    \draw (11.3,-1.4) to [out=120, in=-46] (10.5,-0.46);
    \draw (10.5,-0.46) to [out=45, in=-100] (11.2,0.7);

    \draw [thick] (7.6, 1.5) -- (10.5,2.2);
    \draw [thick] (7.6,-1.2) -- (10.5,-0.46);
    \draw (7.6, 3.5) -- (10.5,4.2);
    \draw (8.4, -2.1) -- (11.3,-1.4);
    \draw (6.8, -2.1) -- (9.7,-1.4);
    \draw (8.22, 3.1) -- (11.12,3.8);
    \draw (8.4, 2.4) -- (11.3,3.1);
    \draw [dashed] (6.9, 0.6) -- (9.8,1.3);
    \draw [dashed] (7, -0.4) -- (9.9,0.3);
    \draw [dashed] (7.3, -1.5) -- (10.2,-0.8);

    \draw (10.97,4) to [out=-5, in=135] (13.2,3.2);
    \draw (13.2,3.2) to [out=-45, in=70] (14,-2.1);
    \draw (14, -2.1) to [out=250, in=-5] (11.9, -2.9);
    \draw (11.9, -2.9) to [out=175, in=-40] (9.4,-1.83);

    \draw (10.5,4.2) to [out=5, in=125] (14,3);
    \draw (14,3) to [out=-55, in=80] (14.46, 0.5);
    \draw (14.46, 0.5) to [out=260, in=40] (13.5,-2.2);
    \draw (13.5, -2.2) to [out=220, in=-10] (11.8,-2.4);
    \draw (11.8, -2.4) to [out=170, in=-40] (10.1,-1.68);
    \draw [thick] (9.4,-1.83) to [out=90, in=180] (9.8,-1.48) to [out=0, in= 90] (10.1,-1.68);

    \draw [thin] (9.4, -1.8) -- (8.9,-1.3);
    \draw [thin] (10.1, -1.7) -- (9.5,-1.2);
    \draw [thin] (8.9, -1.3) to [out=130, in=130] (9.5, -1.2);
    \draw [thick] (9.1, -1.5) to [out=90,in=90] (9.7, -1.4);

    \draw (10,-2.3) to [out=130, in= 130] (10.7, -2.09);
    \draw (10.8,-2.7) to [out=130, in=130] (11.2,-2.28);
    \draw (12.6, -2.5) to [out=0, in=0] (12.6, -2.92);
    \draw (13.35, -2.3) to [out=45, in=45] (13.75, -2.5);
    \draw (11.9, 3.9) to [out=20, in=0] (12,4.19);
    \draw (12.84, 3.5) to [out=-30, in=-30] (13.14, 3.8);
    \draw (13.64, 2.7) to [out=-50, in=-50] (14.07, 2.9);

    \node at (9.1,2) {b};
    \node at (9.1,-0.7) {a};
    \node at (10.1,-1.4) {c};
    \node at (8.8,-1.8) {d};
    \draw [->] (8.95,-1.8) -- (9.3,-1.4);
    \node at (7.4,-1.94) {$\ast$};

    \node at (-2.9,-4.2) {set $\beta_1 = c$};
    \draw (1.7,-4.2) ellipse (2.5 and 0.75);
    \node at (-0.7,-4) {$\ast$};
    \node at (1.7,-4.2) {\tiny $F_0$ empty self-intersection};
    \node at (-1.1,-3.8) {\small new};
    \node at (4.8,-4.2) {$=$};
    \draw [dashed] (5.3,-4) -- (6.7,-3.4);
    \draw [dashed] (6.3,-4.7) -- (7.5,-4.2);
    \draw (6.6,-4.7) -- (7.35,-4.4);
    \draw (6.6, -4.7) to [out=90, in=180] (7,-4.1) to [out=0, in=90] (7.35, -4.4);
    \draw (6.6,-4.7) to [out=145, in=-60] (5.6, -3.7);
    \draw (7.2,-4.15) to [out=145, in=-20] (6.2, -3.4);
    \draw (5.6, -3.7) to [out=120, in=160] (6.2, -3.4);
    \node at (7,-4.53) {$\ast$};
    \node at (6.3,-4) {\small $F_0$};
    \node at (8.1,-3.4) {\tiny over crossings};
    \node at (5.7,-4.7) {\small $I_0$};
    \draw [->] (5.7,-4.5) -- (5.7,-4);
    \draw [->] (8,-3.6) -- (6.68,-3.55);
    \draw [->] (8, -3.6) -- (7.4, -4.1);

    \node at (-2.4,-5.3) {\large $F$};
    \draw [->] (-2.1,-5.1) -- (-0.7,-4.8);
    \draw [->] (-2.1,-5.4) -- (-0.7,-5.6);
    \draw (1.7,-6.1) ellipse (2.5 and 0.75);
    \draw (-0.2,-5.6) -- (-0.2,-6.58);
    \draw (0.6,-5.42) -- (0.6,-6.78);
    \draw (2.3,-5.38) -- (2.3,-6.82);
    \draw (3.2,-5.5) -- (3.2,-6.7);
    \node at (-0.4,-6.1) {$a^\prime$};
    \node at (0.4,-6.1) {$b^\prime$};
    \node at (2.5,-6.1) {$b$};
    \node at (3.4,-6.1) {$a$};
    \node at (1.7,-6.1) {\small $F_1$};
    \node at (6.1,-6) {$=$};

    \draw (7.8,-7.5) to (8.5,-6.9);
    \draw (8.5,-6.9) to [out=45, in=-90] (8.9,-6.3);
    \draw (8.9,-6.3) to [out=90, in=-45] (8.5,-5.7) to [out=135, in=-90] (8.1,-5.1);
    \draw (8.1, -5.1) to [out=90, in=180] (8.5, -4.7);
    \draw (8.5, -4.7) to [out=0, in=90] (8.9, -5.1);
    \draw (8.9, -5.1) to [out=-90, in=45] (8.5, -5.7) to [out=225, in=90] (8.1, -6.3);
    \draw (8.1, -6.3) to [out=-90, in=135] (8.5, -6.9) to (9.2,-7.5);

    \draw (9.6,-7.1) to (10.3,-6.5);
    \draw (10.3,-6.5) to [out=45, in=-90] (10.7,-5.9);
    \draw (10.7,-5.9) to [out=90, in=-45] (10.3,-5.3) to [out=135, in=-90] (9.9,-4.7);
    \draw (9.9,-4.7) to [out=90, in=180] (10.3,-4.3);
    \draw (10.3,-4.3) to [out=0, in=90] (10.7,-4.7);
    \draw (10.7,-4.7) to [out=-90, in=45] (10.3,-5.3) to [out=225, in=90] (9.9,-5.9);
    \draw (9.9,-5.9) to [out=-90, in=135] (10.3,-6.5) to (11,-7.1);

    \draw (8.5, -4.7) -- (10.3, -4.3);
    \draw (8.5, -5.7) -- (10.3, -5.3);
    \draw (8.5, -6.9) -- (10.3, -6.5);
    \draw (7.8,-7.5) -- (9.6, -7.1);
    \draw (9.2, -7.5) -- (11, -7.1);
    \node at (7.8,-7.5) {$\ast$};

    \draw (10.64, -4.5) arc (100:-90:0.9 and 1.6);
    \draw (10.5, -4.35) to [out=20, in=65] (11.67,-5.7);
    \draw (11.67,-5.7) to [out=245, in=85] (11.5,-6.1);
    \draw (11.5,-6.1) to [out=265, in=0] (10.9,-7.4);
    \draw (10.9, -7.4) to [out=160, in=90] (10.5, -7.5);
    \draw (10.5, -7.5) to [out=-90, in=200] (10.8, -7.67);
    \node at (9.4,-5) {$F_1$};

    \node at (-2.4,-9.5) {\large $F_1$};
    \draw [->] (-2.1,-9.4) -- (-0.7,-9);
    \draw [->] (-2.1,-9.6) -- (-0.7,-9.8);
    \node at (-0.6,-5.8) {$\ast$};
    \node at (-1.1,-5.8) {\small old};
    \draw (1.7,-8.6) ellipse (2 and 0.75);
    \draw (1.1,-7.9) -- (1.1,-9.3);
    \node at (-0.3,-8.6) {$\ast$};
    \node at (-0.6,-8.3) {\small new};
    \draw (2.2,-7.9) -- (2.2,-9.3);
    \node at (0.8,-8.6) {$b^\prime$};
    \node at (2.5,-8.6) {$b$};
    \draw (1.7, -10.4) ellipse (2.5 and 0.75);
    \node at (1.7,-10.4) {\tiny $F_{11}$ empty self-intersection};
    \node at (-0.8,-10.4) {$\ast$};
    \node at (-1.3,-10.6) {\small old};
    \draw [->] (4.4,-8.5) -- (6,-8.5);

    \draw (7.5,-10.1) -- (6.8,-9.5);
    \draw (6.8,-9.5) to [out=130, in=-90] (6.7,-9.2);
    \draw (6.7,-9.2) to [out=90, in=225] (7.1,-8.7) to [out=45, in=-90] (7.4,-8.2);
    \draw (7.4,-8.2) to [out=90, in=0] (7.1,-7.9) to [out=180, in=90] (6.7,-8.2);
    \draw (6.7,-8.2) to [out=-90, in=135] (7.1,-8.7) to [out=-45, in=90] (7.4,-9.1);
    \draw (7.4,-9.1) to [out=-90, in=90] (7.2,-9.5);
    \draw (7.2,-9.5) to [out=-70, in=130] (7.6,-9.9);
    \draw (8.05,-10.25) -- (7.6,-9.9);
    \draw (7.1,-10) to [out=130, in=0] (6.7,-9.7) to [out=180, in=50] (6.3,-10);

    \draw (10.15,-9.7) -- (9.45,-9.1);
    \draw (9.45,-9.1) to [out=130, in=-90] (9.35,-8.8);
    \draw (9.35,-8.8) to [out=90, in=225] (9.75,-8.3) to [out=45, in=-90] (10.05,-7.8);
    \draw (10.05,-7.8) to [out=90, in=0] (9.75,-7.5) to [out=180, in=90] (9.35,-7.8);
    \draw (9.35,-7.8) to [out=-90, in=135] (9.75,-8.3) to [out=-45, in=90] (10.05,-8.7);
    \draw (10.05,-8.7) to [out=-90, in=90] (9.85,-9.1);
    \draw (9.85,-9.1) to [out=-70, in=130] (10.25,-9.5);
    \draw (10.7,-9.85) -- (10.25,-9.5);
    \draw [dashed] (9.75,-9.6) to [out=130, in=0] (9.35,-9.3) to [out=180, in=50] (8.95,-9.6);

    \draw (8.05, -10.25) -- (10.7, -9.85);
    \draw (7.1, -8.7) -- (9.75,-8.3);
    \node at (8.4,-8.3) {$b$};
    \draw [dashed] (7.5, -10.1) -- (10.15, -9.7);
    \draw [dashed] (7.1, -10) -- (9.75, -9.6);
    \draw [dashed] (6.3, -10) -- (9, -9.6);
    \node at (6.3,-10) {$\ast$};
    \draw [thin] (7,-9.5) circle (0.4);
    \draw [thin] (9.6,-9.1) circle (0.4);

    \draw [->] (4.5, -10.3) to [out=30, in=180] (5.9, -10);
    \node at (7,-10.4) {$I_{11}$};
    \draw [->] (7, -10.2) to [out=100, in=20] (6.5, -10.05);
    \node at (7.7,-9.2) {\small $F_{10}$};
    \draw (11.7, -9.3) circle (0.3);
    \draw (14.2,-9.3) circle (0.3);
    \node at (12.95,-9.2) {\small over-};
    \node at (12.95,-9.5) {\small crossings};
    \draw (11.5,-9.5) -- (11.9,-9.1);
    \draw [->] (12.4, -9.3) -- (12.1,-9.3);
    \draw (11.53,-9.13) -- (11.65,-9.25);
    \draw (11.75,-9.35) -- (11.87,-9.47);
    \draw (14,-9.1) -- (14.4,-9.5);
    \draw [->] (13.6,-9.3) -- (14.05,-9.3);
    \draw (14,-9.47) -- (14.15,-9.35);
    \draw (14.23,-9.27) -- (14.37,-9.15);
    \node at (12.9,-10.13) {\small (middle sheet compressed};
    \node at (12.9,-10.6) {\small toward int($M$))};

    \node at (-2.3,-13.1) {$F_{10}$};
    \draw [->] (-1.9,-12.9) -- (-0.7,-12.5);
    \draw [->] (-1.9,-13.2) -- (-0.7,-13.5);
    \draw (1.7, -12.1) ellipse (2.5 and 0.75);
    \node at (1.7,-12.1) {\tiny $F_{100}$ empty self-intersection};
    \node at (-0.8,-12.1) {$\ast$};
    \node at (-1.1,-11.8) {\small new};
    \draw (1.7,-13.9) ellipse (2.5 and 0.75);
    \node at (1.7,-13.9) {\tiny $F_{101}$ empty self-intersection};
    \node at (-0.8,-13.9) {$\ast$};
    \node at (-1.2,-14.1) {\small old};
    \node at (5.3,-12.1) {$=$};
    \node at (5.3,-13.9) {$=$};

    \draw (7.9,-14.4) -- (7.4,-13.6);
    \draw (7.4,-13.6) to [out=120, in=-90] (7.3,-13.2);
    \draw (7.3,-13.2) to [out=90, in=145] (7.9,-12.4);
    \draw (7.9,-12.4) to [out=-45, in=40] (7.8,-13.3);
    \draw (7.8,-13.3) to [out=220, in=120] (8.2,-14.3);

    \draw (8.1,-12.35) to [out=-45, in=40] (8,-13.25);
    \draw (8,-13.25) to [out=220, in=120] (8.4,-14.29);
    \draw (8.1, -12.35) to [out=135, in=-90] (7.5, -11.7);
    \draw (7.5, -11.7) to [out=90, in=180] (8,-11.1);
    \draw (8, -11.1) to [out=0, in=90] (8.5, -11.7) to [out=-90, in=45] (8.3,-12);
    \draw (8.3, -12) to [out=225, in=135] (8.3,-12.4);
    \draw (8.3, -12.4) to [out=-45, in=40] (8.2, -13.3);
    \draw (8.2, -13.3) to [out=220, in=120] (8.5,-14.26);
    \draw [thin] (8, -12.1) circle (0.4);

    \draw (10,-13.94) -- (9.5,-13.2);
    \draw (9.5,-13.2) to [out=120, in=-90] (9.4,-12.8);
    \draw (9.4,-12.8) to [out=90, in=145] (10,-12);
    \draw (10,-12) to [out=-45, in=40] (9.9,-12.9);
    \draw (9.9,-12.9) to [out=220, in=120] (10.3,-13.86);
    \draw (10.2,-11.95) to [out=-45, in=40] (10.1,-12.85);
    \draw (10.1,-12.85) to [out=220, in=120] (10.5,-13.83);
    \draw (10.2,-11.95) to [out=135, in=-90] (9.6,-11.3);
    \draw (9.6,-11.3) to [out=90, in=180] (10.1,-10.7);
    \draw (10.1,-10.7) to [out=0, in=90] (10.6,-11.3) to [out=-90, in=45] (10.4,-11.6);
    \draw (10.4,-11.6) to [out=225, in=135] (10.4,-12);
    \draw (10.4,-12) to [out=-45, in=40] (10.3,-12.9);
    \draw (10.3,-12.9) to [out=220, in=120] (10.6,-13.8);
    \draw [thin] (10.1,-11.7) circle (0.4);

    \draw [thin] (7.9, -14.4) -- (10.6, -13.8);
    \draw [dashed] (8, -11.1) -- (10.1, -10.7);
    \draw [dashed] (8.5, -11.7) -- (10.6, -11.3);
    \draw [dashed] (8.2, -12.2) -- (10.3, -11.8);
    \node at (6.9,-12.7) {$F_{101}$};
    \node at (7.9,-14.4) {$\ast$};
    \node at (8.6,-14.3) {$\ast$};
    \node at (7.9,-14.65) {\tiny 101};
    \node at (8.9,-14.55) {\tiny 100};
    \node at (8.3,-15) {\tiny $I_{101}$};
    \draw [->] (8.3,-14.85) -- (8.2,-14.4);
    \node at (8.8,-15) {\tiny $I_{100}$};
    \draw [->] (8.8,-14.85) -- (8.4,-14.4);
    \node at (8.9,-11.3) {$F_{100}$};
    \node at (5.3,-15.6) {$I = I_{11} \circ I_{101} \circ I_{100} \circ I_0$};
    \node at (1.7,-8.6) {\small $F_{10}$};
    \end{tikzpicture}
    \caption{}
  \end{figure}

  We now finish the proof of Theorem 0.

  At the top of the tower the collapse of $((l_n (\pa D))^\pr, \ast)$ is guided by an initial sequence of properly embedded 2-cells $\{E_1, ..., E_j\}$, where the base point of each $E_1, ..., E_j$ is taken to be the self intersection point of its boundary, and with the original $\ast \subset E_j$ serving as the base point for the final $E_j$.

  Now consider each of these cells $E_i$, $1 \leq i \leq j$, and its singular image $F_i$ under the 2-fold cover $\pi_{n-1}: \widetilde{N}_{n-1} \ra N_{n-1}$. Again we may remove double loops and only need consider double arcs. As above, \underln{dissect} $F_i$ into a sequence of great disks $\{F_{i, 1}, ..., F_{i, k_i}\}$ and replace each $E_i$ with that ordered sublist to obtain a great sequence:

  $$\{F_{1,1}, ..., F_{1, k_1}, ..., F_{j, 1}, ..., F_{j, k_j}\}$$

  \noindent
  which guides the collapse of $(l_{n-1} (\pa D))^\pr$ to $\ast$ (and defines the crossing resolution indicated by the "prime").

  Now proceed all the way down the tower. At height $m$ by induction assume a great sequence $\{G_1, ..., G_p\}$ for $(l_m (\pa D))^\pr$. Project each $G_i$, $1 \leq i \leq p$, by $\pi_{m-1}$ to $N_{m-1}$ and construct, as above, a great sequence for each $(\pi_{m-1} (\pa G_i))$. The great sequence for $(l_{m-1} (\pa D))$ is obtained by concatenating these $p$ great sequences. Finally, setting $m = 1$ we obtain a great sequence guiding the null isotopy which establishes Theorem 0.
\end{proof}

\underln{Note}: There is a potentially exponential efficiency in describing a null isotopy as a long composition, as we have done in the proof of Theorem 0, compared with attempting a direct description of the bounding disk. The Fox-Artin-like unknot in Figure 7 makes this clear: any bounding disk must pass $2^k$ times over the rightmost protuberance. This motivated our proof strategy.

\begin{figure}[H]
  \centering
  \begin{tikzpicture}[scale = 1.4]
    \draw (0,0) to [out=-5, in=190] (3.3,0.1);
    \draw (0,0) to [out=175, in=270] (-0.8,1.1);
    \draw (-0.8,1.1) to [out=90, in=160] (0.9,2.5);
    \draw (0.25,1.4) to [out=70, in=180] (1.7, 2.6) to [out=0, in=160] (2,2.55);
    \draw (0.6, 1.5) to [out=70, in=200] (1.2,2.2) to [out=20, in=180] (2.2,1.9);
    \draw (1.9,2.2) to [out=70, in=170] (2.8,2.6);
    \draw (2.2, 1.9) to [out=0, in=270] (2.43, 2.1) to [out=90, in=-30] (2.35, 2.35);
    \draw (2.1,2.1) to [out=70, in=210] (2.35, 2.35) to [out=30, in=180] (3.1,2.2);
    \draw (3.1, 2.2) to [out=0, in=270] (3.25, 2.4) to [out=90, in=-30] (3.2, 2.5);

    \draw (2.8, 2.45) to [out=70, in=200] (3.1, 2.7);
    \draw (3, 2.4) to [out=70, in=180] (3.2, 2.57);

    \draw (1.2,2.1) to [out=-60, in=30] (0.7,1.4);
    \draw (0.7, 1.4) to [out=210, in=70] (-0.3,1);
    \draw (-0.3,1) to [out=250, in=205] (0,0.7);
    \draw (0,0.7) to [out=15, in=250] (0.2,1.1);
    \draw (0.5, 1.2) to [out=250, in=150] (0.6,0.9) to [out=-25, in=250] (1.8,1.9);
    \draw (2, 1.8) to [out=250, in=150] (1.9,0.9) to [out=-25, in=250] (2.75, 2.15);

    \node at (3.35,2.75) {...};
    \node at (3.4,2.55) {...};
    \draw (3.7,2.65) to [out=70, in=160] (4.6,3.05);

    \draw (3.9,2.55) to [out=70, in=210] (4.15,2.8) to [out=30, in=180] (4.9,2.55);
    \draw (4.9,2.55) to [out=0, in=270] (5.1,2.8) to [out=90, in=-30] (5,2.95);

    \draw (3.58,2.55) to [out=-10, in=180] (4.08,2.35);
    \draw (4.08,2.35) to [out=0, in=270] (4.28,2.6) to [out=90, in=-30] (4.18,2.75);
    \draw (3.55, 2.8) to [out=35, in= 170] (3.75, 2.85);
    \draw (2.9,2.1) to [out=250, in=150] (2.9,1.1) to [out=-25, in=250] (3.3,1.4);
    \node [rotate=70] at (3.4,1.6) {...};
    \draw (3.5, 1.8) to [out=70, in=260] (3.65, 2.4);

    \draw (4.8,2.8) to [out=70, in=210] (5.05,3.1) to [out=30, in=180] (5.8,2.85);
    \draw (5.8,2.85) to [out=0, in=270] (6,3.05) to [out=90, in=-30] (5.9,3.25);
    \draw (4.6, 2.9) to [out=70, in=150] (5.9, 3.25);
    \draw (3.8,2.3) to [out=250, in=140] (3.8,1.3) to [out=-35, in=250] (4.5,2.6);

    \draw (3.9,0.2) to [out=20, in=250] (4.7, 2.5);
    \draw  (6.2,2.8) arc (-90:90:0.075 and 0.25);
    \draw  (6.2,2.8) arc (270:90:0.075 and 0.25);
    \draw  (6.3,2.8) arc (-90:90:0.075 and 0.25);
    \draw [dashed] (6.3,2.8) arc (270:90:0.075 and 0.25);
    \node at (6.6,3.05) {...};
    \draw  (6.9,2.8) arc (-90:90:0.075 and 0.25);
    \draw  (6.9,2.8) arc (270:90:0.075 and 0.25);
    \draw  (7,2.8) arc (-90:90:0.075 and 0.25);
    \draw [dashed] (7,2.8) arc (270:90:0.075 and 0.25);

    \node at (1.7,2.9) {1};
    \node at (2.6,2.9) {2};
    \node at (4.2,3.4) {$k$-1};
    \node at (5.4,3.6) {$k$};
    \node at (6.6,2.4) {$2^k$ sheets};
    \node at (6.6,1.9) {of spanning disk};
    \node [rotate=15] at (3.6,0.15) {...};
  \end{tikzpicture}
  \caption{}
\end{figure}

\clearpage

\bibliography{library}{}


\end{document}